\RequirePackage{fix-cm}
\documentclass[10pt]{article}
\usepackage{graphicx}

\usepackage{array}
\usepackage{amsmath}
\usepackage{amssymb}
\usepackage{url}
\usepackage{cite}
\usepackage{breakurl}
\usepackage{multirow}
\usepackage{rotating}
\usepackage{indentfirst}
\usepackage{mdwlist}
\usepackage{setspace}

\usepackage{algpseudocode}
\usepackage{algorithm}

\usepackage{mathtools}

\newcommand{\argmax}{\operatornamewithlimits{argmax}}
\newcommand{\pmf}{\operatorname{pmf}}

\begin{document}

\title{A fast numerical method for max-convolution and the application
  to efficient {\tt max-product} inference in Bayesian
  networks}

\author{Oliver Serang \\
	Freie Universit\"at Berlin and\\
        Leibniz-Institute of Freshwater Ecology and Inland Fisheries (IGB)\\
        orserang@uw.edu
}

\date{\bf \noindent Published in the Journal of Computational Biology, 2015 Aug;22(8):770-83. doi:10.1089/cmb.2015.0013}

\maketitle

\begin{abstract}
\noindent Observations depending on sums of random variables are common
throughout many fields; however, no efficient solution is currently
known for performing {\tt max-product} inference on these sums of
general discrete distributions ({\tt max-product} inference can be
used to obtain \emph{maximum a posteriori} estimates). The limiting
step to {\tt max-product} inference is the max-convolution problem
(sometimes presented in log-transformed form and denoted as ``infimal
convolution'', ``min-convolution'', or ``convolution on the tropical
semiring''), for which no $O(k\log(k))$ method is currently
known. Here I present a $O(k\log(k))$ numerical method for estimating
the max-convolution of two nonnegative vectors (\emph{e.g.}, two
probability mass functions), where $k$ is the length of the larger
vector. This numerical max-convolution method is then demonstrated by
performing fast {\tt max-product} inference on a convolution tree, a
data structure for performing fast inference given information on the
sum of $n$ discrete random variables in $O(n k\log(n k)\log(n) )$
steps (where each random variable has an arbitrary prior distribution
on $k$ contiguous possible states). The numerical max-convolution
method can be applied to specialized classes of hidden Markov models
to reduce the runtime of computing the Viterbi path from $n k^2$ to $n
k\log(k)$, and has potential application to the all-pairs shortest
paths problem.
\end{abstract}

\section{Introduction}
\label{intro}
In many fields it is common to have access to information about sums
of random variables and to desire information about those variables
themselves. In mass spectrometry, when two (or more) analytes with
similar mass-to-charge are measured, the intensity of the resulting
peak is a function of the sum of abundances of those analytes (this
problem occurs not only in the mass spectrometry of small molecules,
but also in measuring isotope measurement in elemental and nuclear
mass spectrometry). In transcriptomics, the abundance of a particular
non-unique read (\emph{i.e.}, an RNA sequence that maps to multiple
locations in the transcriptome or genome) provides information about
the sum of the abundances of all transcripts that contain the read
(each transcript weighted by how many copies of the read it
carries). Proteomics has its own version of non-unique reads, shared
peptides which can be found in multiple proteins (not only are shared
peptides the principal source of difficulty in protein
inference\cite{serang:review, serang:efficient, serang:faster}, they
are also responsible for the difficulty evaluating putatative sets of
discovered proteins\cite{serang:recognizing, serang:npci}). In
population genetics, the prior knowledge about population structure
can suggest an expected number of individuals with a particular
genotype, which in turn yields probabilistic information about the
individuals whose aggregate genotypes are expected to produce that sum
(inference is particularly pronounced in polyploids, which increase
the dimensionality of the problem\cite{serang:efficient2}).

In all of these fields, the information on sums of random variables
presents a singular obstacle to computational biology. And regardless
of how infrequently we as scientists directly discuss our current
inability to effectively utilize information about sums of random
variables, the perception of our limited ability to meet the challenge
has become firmly entrenched in our collective unconscious; the silent
agreement on our inability to turn the sausage grinder backwards and
convert the sausage (information about the sum of several random
variables) back into the pigs (information about those random
variables that contributed to the sum) is so well-established that it
not only defines the way that we address these data (\emph{i.e.}, mass
spectra peaks containing multiple analytes, counts of non-unique
reads, \emph{etc.}), but it also causes us to discard data and even
limit research directions we might otherwise consider. For instance,
in mass spectrometry, substantial effort is invested in
chromatography\cite{barnes1992high, james1952gas} and other separation
techniques\cite{pringle2007investigation}, which aim to distinguish
and separate analytes so that they will not be measured by the mass
spectrometer simultaneously (thereby reducing the chances of analytes
resulting in overlapping peaks); the task of decomposing this useful
aggregate information from overlapping peaks back into information
about its contributing parts is eschewed in favor of significant
investments in instrumentation (\emph{e.g.}, more and more advanced
separation technologies and higher-resolution mass
spectrometers\cite{polacco2011discovering}), and still it is common
practice to discard shared isotope peaks and shared peptides even
though they may comprise a large percent of the data and contain
additional information\cite{dost:shared, serang:npci}. Likewise, in
genomics and transcriptomics it is common to simply discard all data
from non-unique reads\cite{lefranccois2009efficient,
  zentner2011integrative}. In burgeoning fields such as metagenomics,
this loss of data-- and subsequent loss of information-- can be even
more pronounced, for instance when all data that map to two or more
species of interest are discarded. In some cases, recovering this lost
information will be the key to making strong conclusions, such as
distinguishing between two closely related species (or bacterial
strains) in a metagenomic mixture.

\subsection*{{\tt Sum-product} inference}
Fast Fourier transform- (FFT-)based convolution can be used to
dramatically improve the efficiency of computing the {\tt sum-product}
addition of two discrete random variables. For three discrete random
variables where $M = L+R$ and where $L \in \{ 0, 1, \ldots k-1 \}$ and
$R \in \{ 0, 1, \ldots k-1 \}$, then the probability mass function
(PMF) of $M$ can be computed via the convolution of the PMFs of $L$
and $R$, denoted $\pmf_L$ and $\pmf_R$ respectively. Note that it is
sufficient to compute $\pmf_M' \propto \pmf_M$, because the result
will be scaled so that its sum is of unity ($\pmf_M = \frac{\pmf_M'}{
  \sum_m \pmf_M'[m]}$):
\begin{eqnarray*}
\pmf_M[m] & \propto & \sum_\ell \sum_r \Pr(L=\ell) \Pr(R=r) \Pr(M=\ell+r) \\
&=& \sum_\ell \Pr(L=\ell) \Pr(R=m-\ell)\\
&=& \pmf_{L} * \pmf_{R} \\
\end{eqnarray*}
where $*$ performs the convolution between the two $k-$length vectors
storing each PMF. 

Whereas naive convolution would compute $\pmf_{M}$ in $O(k\times k)$
steps, FFT exploits the bijection of this convolution to the product
between two polynomials (where the vectors being convolved are the
coefficients of the polynomials being multiplied and the coefficients
of their product forms the vector result); this bijection enables the
use of alternative forms for representing the polynomials (each order
$k-1$ polynomial can be represented through $k$ unique points through
which it passes), which in turn permits elegant divide and conquer
algorithms such as the Cooley-Tukey FFT to compute fast convolution in
$k\log(k)$ steps. Subtraction in the {\tt sum-product} (\emph{i.e.},
computing $L = M-R$) scheme can be performed by first negating $R' =
-R$ (this is done by reversing reversing the vector storing the PMF
$\pmf_{R'}[r] = \pmf_{R}[-r]$), and then adding $L=M+R'$ as before via
FFT convolution $\pmf_{L}' = \pmf_{M} * \pmf_{R'}$. Also, the runtime
constant on FFT convolution algorithms is generally very low, partly
due to the nature of elegance of the algorithms and partly because
implementations have been optimized heavily due to the ubiquity of
convolution in signal processing.

The task of processing information about the sum of $n$ discrete
variables (each with $k$ bins) to retrieve information on the
individual variables can be performed naively in $O(k^n)$ steps by
simply enumerating the exponentially many possible outcomes; however,
such brute-force techniques are wildly inefficient when either $n$ or
$k$ become large. Fortunately, recent work proposes methods to
decompose larger problems (\emph{e.g.}, into multiple sums and
differences of pairs of discrete variables of the form $M=L+R$ and
$L=M-R$, very fast inference can be achieved. This has been derived
for binary variables ($k=2$) in
$n\log(n)\log(n)$\cite{tarlow2012fast}, and was independently
discovered for arbitrary discrete distributions (\emph{i.e.}, where
$k>2$) and to multidimensional distributions (via matrix convolution,
which can be decomposed into one-dimensional convolutions by the
row-column algorithm) using the probabilistic convolution
tree\cite{serang2014probabilistic} ({\bf
  algorithm~\ref{algorithm:convolution-tree}}). In the general case,
distributions on all individual variables conditional on information
about the sum can be computed in $O(n k \log(n k) \log(n))$ steps
(whenever $k\log(k)$ fast convolution is available):
\begin{eqnarray*}
convolutionTreeCost(n,k) & = & \sum_{u=1}^{log(n)} \frac{n}{2^u} ~ convolutionCost(k 2^u) \\
& = & \sum_{u=1}^{log(n)} \frac{n}{2^u} ~ k 2^u \log(k 2^u) \\
& = & n k \sum_{u=1}^{log(n)} \log(k) + \log(2^u)\\
& = & n k \sum_{u=1}^{log(n)} \log(k) + u\\
& = & n k  \left[ \sum_{u=1}^{log(n)} \log(k) \right] + \left[ \sum_{u=1}^{log(n)} u \right]\\
& = & n k  \left[ \log(n) \log(k) \right] + \left[ \frac{\log(n)(\log(n)+1)}{2} \right]\\
& \in & O(n k\log(n k)\log(n).\\
\end{eqnarray*}
In practice, this can be significantly faster than the $O(n^2 k^2)$
steps required by dynamic programming when fast convolution is not
available. For instance, when an observed transcript fragment could
originate from $n=256$ species, and where the abundance of each
species is discretized into $k=1024$ bins, then fast convolution makes
inference more than 1800 times faster (the difference between one
algorithm taking 1 second and the other taking 30 minutes), and the
disparity only grows for problems with larger values of $n$ and
$k$. It should be noted that these are only approximate runtimes
calculated from the Big O form; in practice, it is fairly likely that
methods based on fast convolution will be significantly faster,
because of the method's inherent properties and the fact that very
optimized implementations exist.

\begin{algorithm}
  \caption{ {\bf The probabilistic convolution tree algorithm}
    utilizes fast convolution (or fast max-convolution) to efficiently
    turn information on sums of variables back into information on the
    individual variables. The first parameter $(\pmf_{X_1},
    \pmf_{X_2}, \ldots \pmf_{X_n})$ is a collection of $n$
    multidimensional discrete distributions (with same dimension). In
    the case of univariate distributions, they are one-dimensional
    PMFs with $k$ possible outcomes. The second parameter $\pmf_M$ is
    a multidimensional discrete distribution (with same dimension as
    $X_1, X_2, \ldots$) where $M = X_1 + X_2 + \ldots + X_n$. The
    third parameter is a convolution operator (\emph{e.g.}, either
    standard convolution or max-convolution). The algorithm returns a
    pair of values. The first is a collection of likelihood
    distributions given on the information from the sum $M$,
    $(\pmf_{Y_1}, \pmf_{Y_2}, \ldots \pmf_{Y_n})$. The second is the
    prior distribution of $M$ after adding together all
    $X_1+X_2+\ldots X_n$. When $n$ is not an integer power of 2, dummy
    variables whose PMFs of length 1 (\emph{i.e.}, PMFs that are
    Kronecker deltas with 100\% chance of having value 0) should be
    padded on until the next power of 2 is reached (this will allow
    construction of a full binary tree without changing the sum).}

  \label{algorithm:convolution-tree}
  \begin{small}
    \begin{algorithmic}[1]
      \Procedure{probabilisticConvolutionTree}{$(\pmf_{X_1}, \pmf_{X_2}, \ldots \pmf_{X_n})$, $\pmf_M$, $*_{param}$}
      
      \State $forwardTree \gets [ ~ [ \pmf_{X_1}, \pmf_{X_2}, \ldots \pmf_{X_n} ] ~ ] $
      \For{$i=1$ to $len(log_2(n))$ $i += 1$}
      \State $currentLayer \gets forwardTree[-1] $ \Comment{Where $[-1]$ gives the last index as in Python}
      \State $newLayer \gets [ ~ ] $
      \For{$j=1$ to $len(currentLayer)-1$, $j += 2$}

      \State $sum \gets currentLayer[j] *_{param} currentLayer[j+1]$
      \State $newLayer.append( normalized(sum) ) )$ \Comment{Add the PMFs using the $*_{param}$ operator}
      \EndFor
      \State $forwardTree.append(newLayer)$
      \EndFor

      \State $reverseTree \gets [ ~ [ \pmf_M ] ~ ] $
      \For{$i=1$ to $len(log_2(n))$, $i += 1$}
      \State $currentLayer \gets reverseTree[-1] $
      \State $newLayer \gets [ ~ ] $

      \State $forwardTreeLayerToSubtract \gets forwardTree[-i-2] $ \Comment{Python notation for the second to last layer in forwardTree}
      \For{$j=1$ to $len(currentLayer)$, $j += 1$}

      \State $lhs \gets forwardTreeLayerToSubtract[2 j]$
      \State $rhs \gets forwardTreeLayerToSubtract[2 j+1]$

      \State $rhsSubtracted \gets currentLayer[j] ~*_{param}~ (-rhs)$ \Comment{Negation reverses the vector}\\
      \Comment{Then the $*_{param}$ operator performs addition}
      \State $newLhs \gets normalized(rhsSubtracted)$ 
      \State $lhsSutracted \gets currentLayer[j] ~*_{param}~ (-lhs)$
      \State $newRhs \gets normalized(lhsSubtracted)$

      \State $newLhs \gets newLhs.narrowToTheSupportOf(lhs)$
      \State $newRhs \gets newRhs.narrowToTheSupportOf(rhs)$

      \State $newLayer.append( newLhs )$ 
      \State $newLayer.append( newRhs )$ 
      \EndFor
      \State $reverseTree.append(newLayer)$
      \EndFor

      \State $\forall j,~ \pmf_{Y_j} \gets reverseTree[-1][j]$
      \State $\pmf_{Z} \gets forwardTree[-1][0]$
      \State \Return $(~(\pmf_{Y_1}, \pmf_{Y_2}, \ldots \pmf_{Y_n}), \pmf_{Z}~)$ \Comment{First return value: collection of likelihood distributions}
      \\ \Comment{Second return value: prior distribution on $M$}
      \EndProcedure
    \end{algorithmic}
  \end{small}
\end{algorithm}

\subsection*{{\tt Max-product} inference}
Qualitatively, {\tt max-product} inference is a close cousin to {\tt
  sum-product} inference. Where {\tt sum-product} inference considers
each of the exponentially many joint events and allows each to
contribute to the result (in hidden Markov models, this is analogous to
the forward-backward algorithm), {\tt max-product} inference allows
only the highest-quality joint events to contribute (in hidden Markov
models, this defines the Viterbi path). Both inference methods have
complementary advantages and disadvantages: The advantage of {\tt
  sum-product} inference is its democratized equal weighting of all
joint events, the variety of which can provide a rich description of
any high-probability joint events suggested by the data; however, this
can also have disadvantages in that many low-quality joint events
(\emph{i.e.}, those with low joint probability) may shape the result
as much as a small number of high-quality results. Likewise, in {\tt
  sum-product} inference, multiple mutually exclusive joint events can
simultaneously contribute to the result, raising the potential to
erroneously infer implausible conclusions, because both may be
plausible before considering the other. It is because of these
disadvantages in {\tt sum-product} inference that {\tt max-product}
inference is widely used, because it forces the inferences to be
jointly plausible (not simply individually, but as a whole), and
because it drowns out noise from low-quality configurations that can
diffuse and lower the certainty of conclusions in {\tt sum-product}
inference.

Specifically, efficient {\tt max-product} inference on sums of random
variables would be quite useful; in addition to the examples of shared
peptides, non-unique reads, \emph{etc.} found throughout computational
biology (wherein we have information about the sum of variables, but
want to draw conclusions about the variables themselves), more
efficient {\tt max-product} inference would make possible new
inference algorithms on specialized classes of hidden Markov models
(HMMs) where the transition probabilities from the state at index $a$
to the state at index $b$ depend on a function of either $a+b$ or
$a-b$ or $b-a$. Such HMMs have applications in finance and time series
analysis, where the $k$ states at any layer are high-resolution
discretizations of some quantity or price, and where the probability
of a price moving up or down is influenced by the quantity up or down
it moved since the previous time point. In a general HMM with $k$
states and $n$ layers of those states, performing either {\tt
  sum-product} (via the forward-backward algorithm) or {\tt
  max-product} (via the Viterbi algorithm) inference requires $O(n
k^2)$ steps; however, performing {\tt sum-product} inference on the
specialized class of HMMs mentioned above would require only $O(n
k\log(k))$ steps, because each layer can be processed as a two-node
convolution tree. But finding the Viterbi path on such a model in $O(n
k\log(k))$ steps is not currently feasible, because doing so would
require performing max-convolution (where the max of all valid
pairings is chosen rather than the sum over all valid pairings) in
$O(k\log(k))$ steps.

However, despite the promise of {\tt max-product} inference on sums of
random variables, a fast practical solution that utilizes
$O(k\log(k))$ max-convolution (\emph{i.e.}, one with speed roughly
comparable to FFT-based standard convolution) is not yet available for
the general max-convolution problem problem. One special cases for use
only when $k=2$\cite{gupta2007efficient, tarlow2010hop}, can solve the
problem in $n\log(n)$ time by sorting the $n$ variables in descending
order of probability $\Pr(X_1 = 1) \geq \Pr(X_2 = 1) \geq \cdots \geq
\Pr(X_n = 1)$ and then exploiting the property that any case where the
number of ``true'' variables $\sum_j X_j = m$ must prefer the first
$m$ variables in the sorted order. This method understandably fails
when $k>2$ because there is no guarantee of an ordering that will
satisfy all dimensions (when $k=2$, increasing the probability of
$\Pr(X_j = 1)$ has a useful effect of decreasing the probability of
$\Pr(X_j = 0)$ by the same amount in order to preserve the unitary
value of the sum). When $k > 2$, a similar idea to the sorting
approach can be used, but not without approximation or some method for
exploring or optimizing the exponential space of joint
events\cite{serang:efficient2}.

Adapting the probabilistic convolution tree algorithm from {\bf
  algorithm~\ref{algorithm:convolution-tree}} (by simply replacing all uses
of $*$ with $*_{\max}$ when adding pairs of variables) achieves only
an $O(n^2 k^2)$ runtime for {\tt max-product} inference because
additions and subtractions between individual pairs of random
variables will require $O(k_1 \times k_2)$ time without a faster
algorithm for max-convolution. It is tempting to try to derive an FFT
equivalent to max-convolution, a subtle difference makes this
challenging: Where standard convolution uses the operations
$(+,\times)$ on real-valued numbers (a ``ring''), max-convolution
employs $(\max,\times)$ (alternatively applying a log-transformation
to the probabilities being convolved will negate them and thus change
the problem to the equivalent $(\min,+)$ operations, called
min-convolution, infimal convolution, or convolution on the ``tropical
semiring''). Regardless of which form is used, the employment of the
$\max$ (or $\min$ in the min-convolution case) downgrades the
operation from a ring to a ``semiring'' because the $\max$ and $\min$
operations have no inverse. Thus, the {\tt max-product} addition of
two discrete random variables takes a different form, which is no
longer bijective to polynomial multiplication:
\begin{eqnarray*}
\pmf_M[m] & \propto & \max_\ell \max_r \Pr(L=\ell) \Pr(R=r) \Pr(M=\ell+r) \\
&=& \max_\ell \Pr(L=\ell) \Pr(R=m-\ell)\\
&=& \pmf_L *_{\max} \pmf_R \\
\end{eqnarray*}
where $*_{\max}$ is the max-convolution operator. The loss of the
bijective polynomial representation prevents the exploitation of the
Lagrange form of polynomials, and thus there is no known $k\log(k)$
algorithm for performing max-convolution. 

Excluding such highly specialized methods as the rank method of
Babai\cite{babai2009computing} (which achieves a runtime of $O(k_1 \sqrt{k_2
  \log(k_2)})$ time but only when the vector of length $k_2$ contains
elements with value $0$ or $\infty$), the two most sophisticated
max-convolution algorithms applicable to probabilistic inference are
from Bussieck \emph{et al.}\cite{bussieck:fast} and Bremner \emph{et
  al.}\cite{bremner:necklaces}.

The method from Bussieck has a $O(k^2)$ runtime in the worst case, but
under a certain distribution values in the two vectors being
convolved, the authors demonstrate an expected runtime of
$O(k\log(k))$\cite{bussieck:fast}. The approach works by starting the
result $\forall m,~\pmf_M'[m] \gets -\infty$ and then proceeds by sorting the
two vectors being convolved ($L$ and $R$) in descending order. Their
method then proceeds through both lists head-first to generate the
first $k\log(k)$ sorted terms of $\forall \ell, r~ \pmf_L[\ell]
\pmf_R[r]$. Each of these terms is used to update the appropriate
index in the result vector $\pmf_M'[\ell+r] \gets \max(\pmf_M'[\ell+r], \pmf_L[\ell]
\pmf_R[r] )$. Thus far, the algorithm is $\in O(k\log(k))$, but there may
be indices of $\pmf_M'$ that have not yet been set (they are still equal to
$-\infty$); each of these must be computed, and each such direct
computation takes $O(k)$ time (if there are $\Omega(k)$ such unset
indices, then the overall runtime becomes $O(k^2)$). Despite the
significant achievement posed by the construction of this algorithm,
the authors suggest that the runtime constant is quite high due to the
overhead of the sophisticated algorithms used to sort the largest
$k\log(k)$ values while neither sorting nor even \emph{generating} all
$k^2$ values; they suggest that their result is of mostly theoretical
import, and suggest using other methods in practice.

The method of Bremner \emph{et al.} (which was subsequently extended
by Williams\cite{williams:faster}) draws a relationship between
min-convolution and the necklace alignment problem, wherein two two
collections of beads, each on its own circular string, are rotated to
optimally align\cite{bremner:necklaces}. Their method is the most
sophisticated in existence, and consists of a highly complicated
exploitation of similarity to the all-pairs shortest paths problem to
achieve a method with a subquadratic worst-case runtime of $O(k^2
\frac{{(\log(\log(k)))}^3}{\log(k)\log(k)})$ for each max-convolution
(the runtime of the Bremner \emph{et al.} method can also be improved
using a more recent method for solving the all-pairs shortest paths
problem, decreasing the runtime to
$\frac{k^2}{2^{\Omega(\sqrt{\log(k)})}}$\cite{williams:faster}). Even
if it were possble to be implemented with runtime constant as
optimized as FFT, the cost of using a max-convolution tree to solve
the previously mentioned metagenomic {\tt max-product} inference
problem on $n=256$ variables where each has $k=1024$ states would be
over 166 times slower than the cost of solving an equally sized {\tt
  sum-product} problem with FFT convolution (the number of steps
required was calculated numerically to avoid computing a closed form
of the computational cost).

Thus, even with significant mathematical sophistication of these two
state-of-the-art methods, practically efficient {\tt max-product}
inference may be out of reach for even moderately sized problems.

\section{A numerical method for efficiently estimating max-convolution}

Here I will introduce a numerical method for estimating the
max-convolution in $k\log(k)$ time, which can be applied easily using
existing high-performance numerical software libraries. This is
essentially performed by transforming both the inputs and outputs of
the FFT to achieve $p-$norm convolution, which in turn is used as an
approximation for max-convolution via the Chebyshev norm. \newline

For a max-convolution between $\pmf_L$ and $\pmf_R$
\[
  \pmf_M'[m] = \max_\ell~ \pmf_L[\ell] \pmf_R[m-\ell],
\]
at each $m$ value, the shifted products terms can be rewritten as a
simple vector $u^{(m)}$ where elements are defined by
\[u^{(m)}[\ell] = \pmf_L[\ell] \pmf_R[m-\ell].\]
Furthermore, because PMFs consist of nonnegative real values (or
machine-precision representations), then this can be rewritten using
the Chebyshev norm, which computes the maximum absolute value in the
vector $u^{(m)}$. Because $u^{(m)}$ comes from the product of PMF
terms, it is also nonnegative and thus absolute values can be ignored
in both the computation and the result:
\begin{eqnarray*}
\pmf_M'[m] &=& \max_\ell u^{(m)}[\ell] \\
|\pmf_M'[m]| &=& \lim_{p \to \infty} \| u^{(m)} \|_p\\
|\pmf_M'[m]| &=& \lim_{p \to \infty} {\left( \sum_\ell |{u^{(m)}[\ell]}|^p \right)}^\frac{1}{p} \\
\pmf_M'[m] &=& \lim_{p \to \infty} {\left( \sum_\ell {u^{m)}[\ell]}^p \right)}^\frac{1}{p}. \\
\end{eqnarray*}
And then each element of $u^{(m)}[\ell]$ can be expanded back into its
original factors:
\[\pmf_M'[m] = \lim_{p \to \infty} {\left( \sum_\ell {\pmf_L[\ell]}^p {\pmf_R[m-\ell]}^p \right)}^\frac{1}{p}.\]
At this point, a sufficiently large value $p^*$ is used in place of
the limit $p \rightarrow \infty$:
\[\pmf_M'[m] \approx {\left( \sum_\ell {\pmf_L[\ell]}^{p^*} {\pmf_R[m-\ell]}^{p^*} \right)}^\frac{1}{p^*}.\]

At this point, it can be observed that every time elements of the PMFs
$\pmf_L[\ell]$ and $\pmf_R[m-\ell]$ appear, they are raised to the
$p^*$ power; thus, it is possible to change variables and let
\begin{eqnarray*}
  \forall \ell &  v_L[\ell] =& {\pmf_L[\ell]}^{p^*}\\
  \forall r &  v_R[r] =& {\pmf_R[r]}^{p^*},\\
\end{eqnarray*}
yielding
\[\pmf_M'[m] \approx {\left( \sum_\ell v_L[\ell] v_R[m-\ell] \right)}^\frac{1}{p^*}.\]

A similar strategy can be made for $\pmf_M'$; it is possible to
introduce another vector $v_M$ such that every element
$\pmf_M'[m]$ is the result of
\begin{eqnarray*}
  \pmf_M'[m] & \approx & {v_M[m]}^{\frac{1}{p^*}}\\
  v_M[m] & = & \sum_\ell v_L[\ell] v_R[m-\ell]\\
\end{eqnarray*}
And thus it becomes clear that $v_M$ is the result of standard
convolution (\emph{not} a max-convolution) between $v_L$ and
$v_R$. This suggests a numerical algorithm that can make use of
existing FFT convolution libraries to compute $v_M = v_L*v_R$ in
$O(k\log(k))$ steps ({\bf algorithm~\ref{algorithm:max-convolution-original}}).

\begin{algorithm}
  \caption{ {\bf Numerical max-convolution (initial version)}, a
    numerical method to estimate the max-convolution of two PMFs or
    nonnegative vectors. The parameters are two PMFs $\pmf_L$ and
    $\pmf_R$ (by definition nonnegative) and the numerical value $p^*$
    used for computation. The return value is a numerical estimate of
    the max-convolution $\pmf_L *_{\max} \pmf_R$.}

  \label{algorithm:max-convolution-original}
  \begin{small}
    \begin{algorithmic}[1]
      \Procedure{maxConvolution}{$\pmf_L$, $\pmf_R$, $p^*$}
      \State $\forall \ell, ~ v_L[\ell] \gets {\pmf_L[\ell]}^{p^*}$

      \State $\forall r, ~ v_R[r] \gets {\pmf_R[r]}^{p^*}$

      \State $v_M \gets \pmf_L ~*~ \pmf_R$ \Comment{Standard FFT convolution is used here}

      \State $\forall m, ~ \pmf_M'[m] \gets {v_M[m]}^{\frac{1}{p^*}}$
      
      \State \Return $\pmf_M'$ \Comment{The return value is an
        estimate of the max-convolution result} \EndProcedure
    \end{algorithmic}
  \end{small}
\end{algorithm}

\subsection{Reducing underflow}

The main caveat for numerical methods is often the loss of precision
due to underflow when raising small probabilities to the power $p^*$
(in this case, no overflow occurs because the inputs are
probabilities); such losses may not be undone later when raising to
the $\frac{1}{p^*}$ power. One way to limit unnecessary loss of
precision is to recognize that since $\pmf_M'$ will be normalized
after it is estimated, then it can be scaled arbitrarily during
computation without altering the final result. For this reason it can
be beneficial to scale a vector by dividing by its maximum element
before raising it to the $p^*$ or $\frac{1}{p^*}$ power; this will
start the dominant elements close to 1, and thus allow them to lose
little information to underflow. Using this strategy yields a slightly
modified algorithm ({\bf
  algorithm~\ref{algorithm:max-convolution-revised}}).

\begin{algorithm}
  \caption{ {\bf Numerical max-convolution (normalized version)}, a
    numerical method to estimate the max-convolution of two PMFs or
    nonnegative vectors (revised to reduce underflow). The parameters
    are two PMFs $\pmf_L$ and $\pmf_R$ (by definition nonnegative) and
    the numerical value $p^*$ used for computation. The return value
    is a numerical estimate of the max-convolution $\pmf_L *_{\max}
    \pmf_R$.}

  \label{algorithm:max-convolution-revised}
  \begin{small}
    \begin{algorithmic}[1]
      \Procedure{maxConvolutionRevised}{$\pmf_L$, $\pmf_R$, $p^*$}
      
      \State $\ell_{\max} \gets \argmax_\ell \pmf_L[\ell]$
      \State $r_{\max} \gets \argmax_r \pmf_R[r]$

      \State $\forall \ell, ~ v_L[\ell] \gets { \left( \frac{\pmf_L[\ell]}{\pmf_L[ \ell_{\max} ]} \right) }^{p^*}$

      \State $\forall r, ~ v_R[r] \gets { \left( \frac{\pmf_R[r]}{\pmf_R[ r_{\max} ]} \right) }^{p^*}$

      \State $v_M \gets \pmf_L ~*~ \pmf_R$ \Comment{Standard FFT convolution is used here}

      \State $m_{\max} \gets \argmax_m v_M[m]$ \Comment{Correct any small errors in FFT convolution}

      \State $\forall m, ~ \pmf_M'[m] \gets { \left( \frac{v_M[m]}{v_M[m_{\max}]} \right) }^{\frac{1}{p^*}}$
      
      \State \Return $\pmf_L[\ell_{\max}] \pmf_R[r_{\max}] ~
      \pmf_M'$ \Comment{Undo the previously performed scaling} \EndProcedure
    \end{algorithmic}
  \end{small}
\end{algorithm}

Note that, like standard implementations of fast convolution $*$, in
implementing the fast $*_{\max}$ operator it is possible to
automatically choose between a naive implementation or the fast
numerical implementation depending on the size of the problem; on very
small problems (\emph{e.g.} $k=8$), the naive operation will have less
overhead and can be a bit faster (the specific threshold can be chosen
roughly by comparing the expected running time from the fast numerical
method $O( k' \log(k') )$ (where $k'$ is double the next integer power
of two and the log is base two) to the $O(k_1 \times k_2)$); this can
reduce numerical error further.

\section{Results}

I briefly compare the speed and accuracy of the fast numerical
max-convolution method as compared to naive max-convolution. Both
methods are implemented in the Python programming language using
floating point math and the {\tt numpy} package (the fast numerical
max-convolution method is implemented from {\bf
  algorithm~\ref{algorithm:max-convolution-revised}}. 

\subsection{Practical efficiency of fast numerical max-convolution}
\begin{figure}
  \centering
  \includegraphics[width=4in]{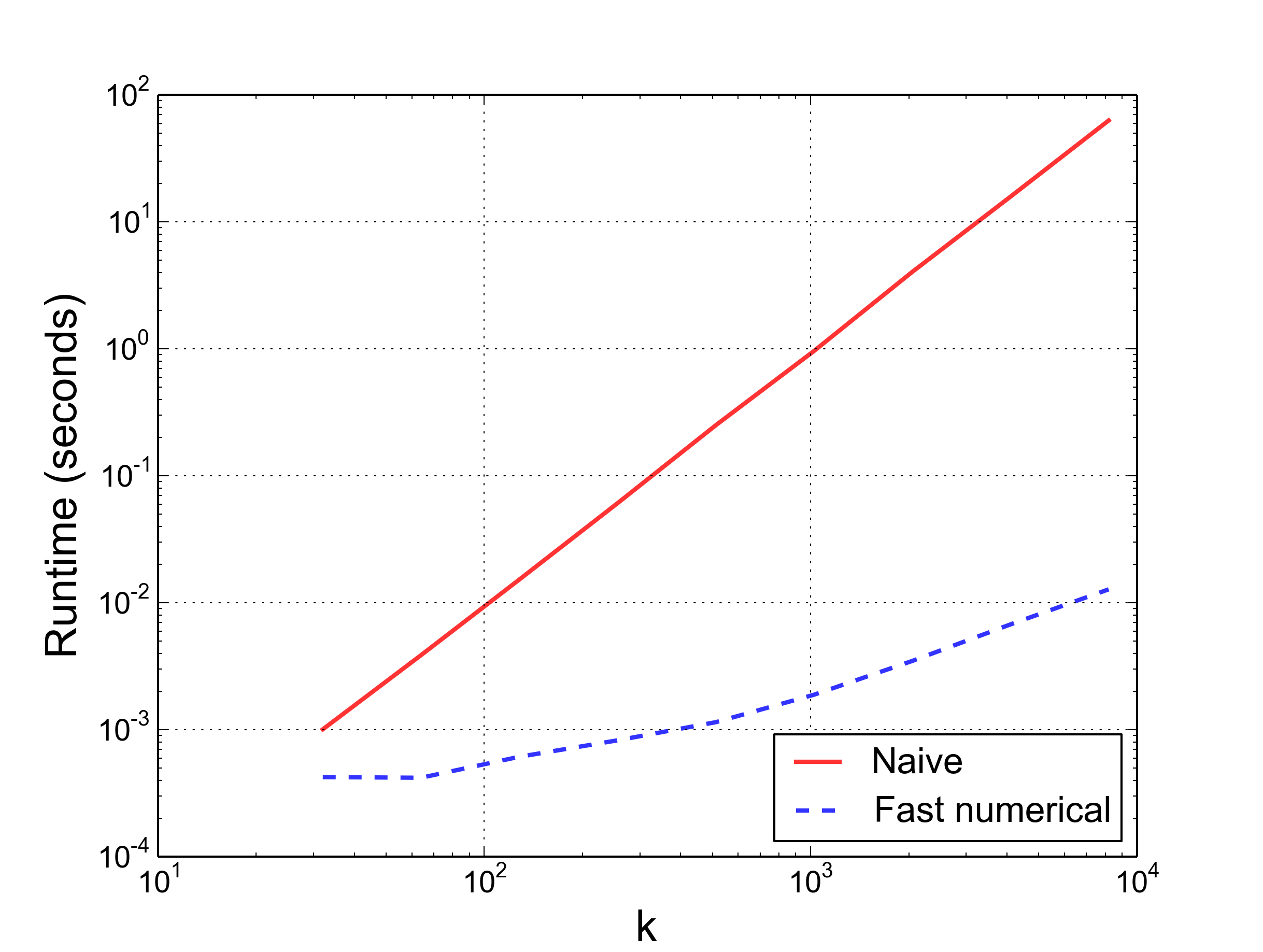} 

  \caption{{\bf Runtime comparison between naive and fast numerical
      convolution.} Note the increasing gap between the two curves,
    which indicates a nonlinear speedup because of the log-scaling of
    both axes (because the runtime of the fast numerical method is
    dominated by FFT calculation).
  \label{figure:max-convolution-speed}}
\end{figure}

The speed of naive max-convolution is compared to the fast numerical
estimate. For each $k \in \{ 32, 64, 128, 256, 512, 1024, 2048, 4096, 8192 \}$, random
pairs of vectors with uniform elements (\emph{i.e.}, each element is
drawn from $uniform(0,1)$). The result of the max-convolution $\pmf_L'
*_{\max} \pmf_R'$ is computed via $O(k^2)$ naive max-convolution and
the fast numerical method. {\bf
  Figure~\ref{figure:max-convolution-speed}} demonstrates an
substantial speedup in practice.

\subsection{Accuracy of fast numerical max-convolution compared to naive max-convolution}
\begin{figure*}
  \centering
  \begin{tabular}{ccc}
    $k=128$ &
    $k=256$ \\
    \includegraphics[width=2in]{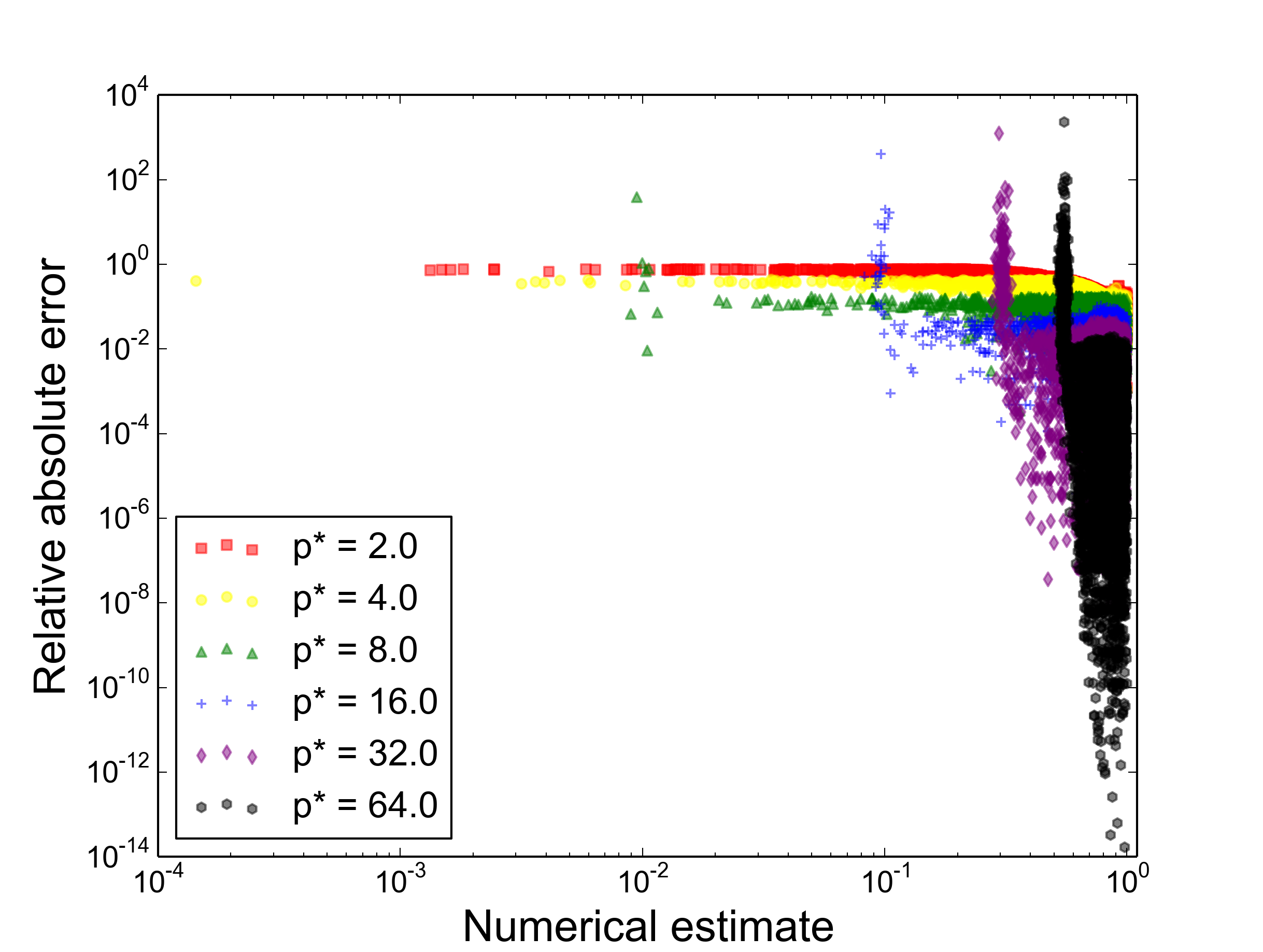} &
    \includegraphics[width=2in]{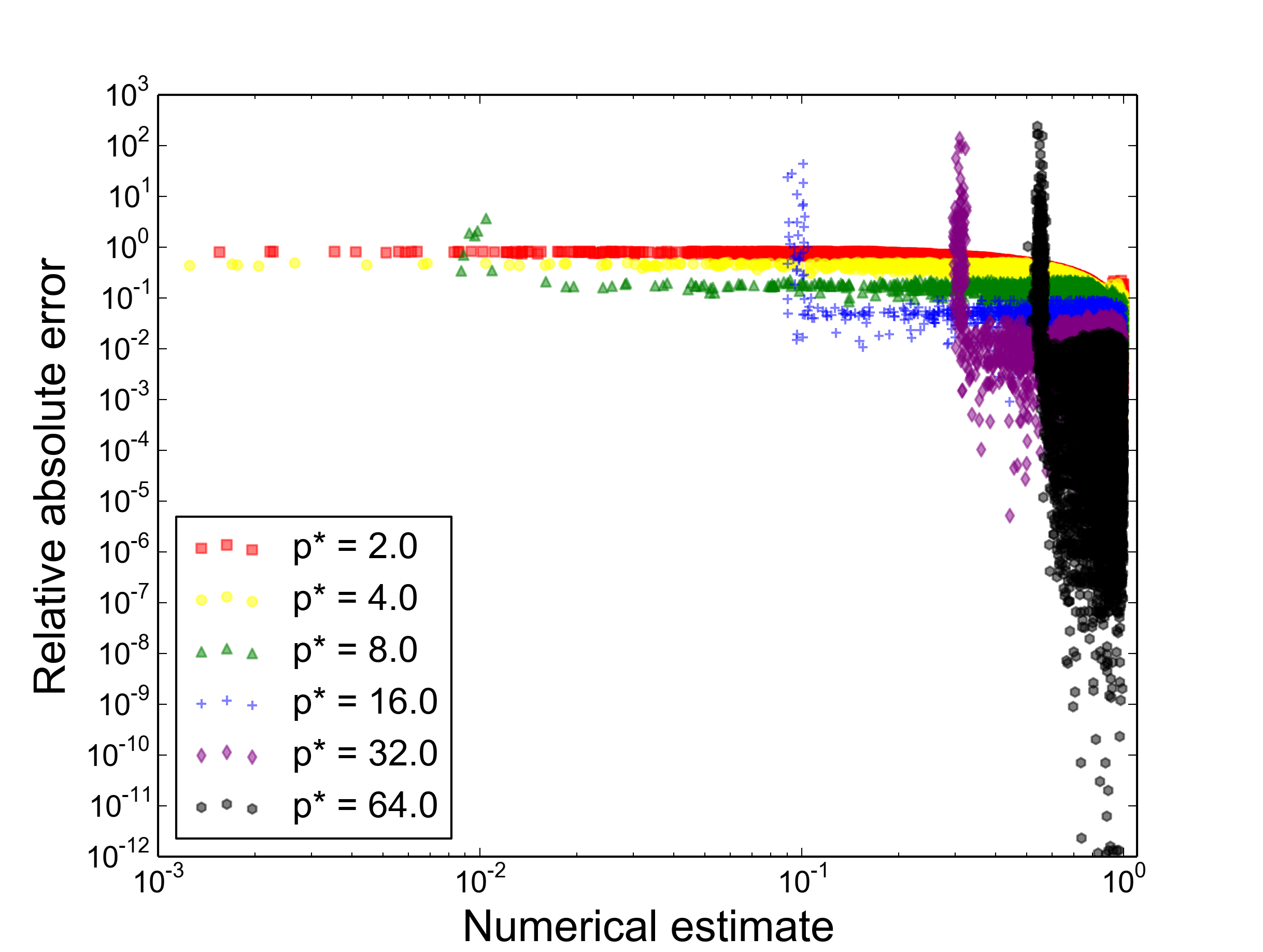} \\
    $k=512$ &
    $k=1024$ \\
    \includegraphics[width=2in]{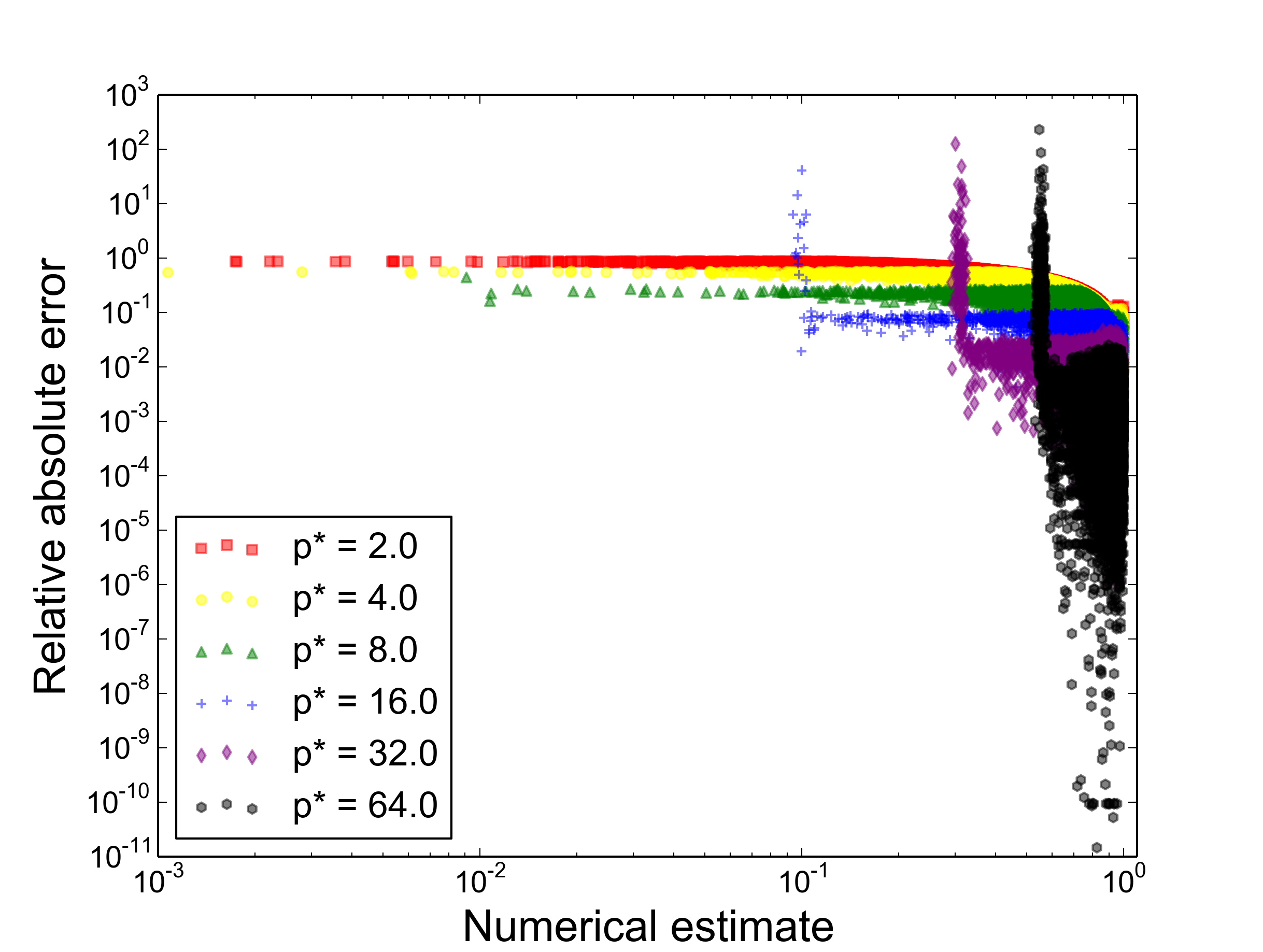} &
    \includegraphics[width=2in]{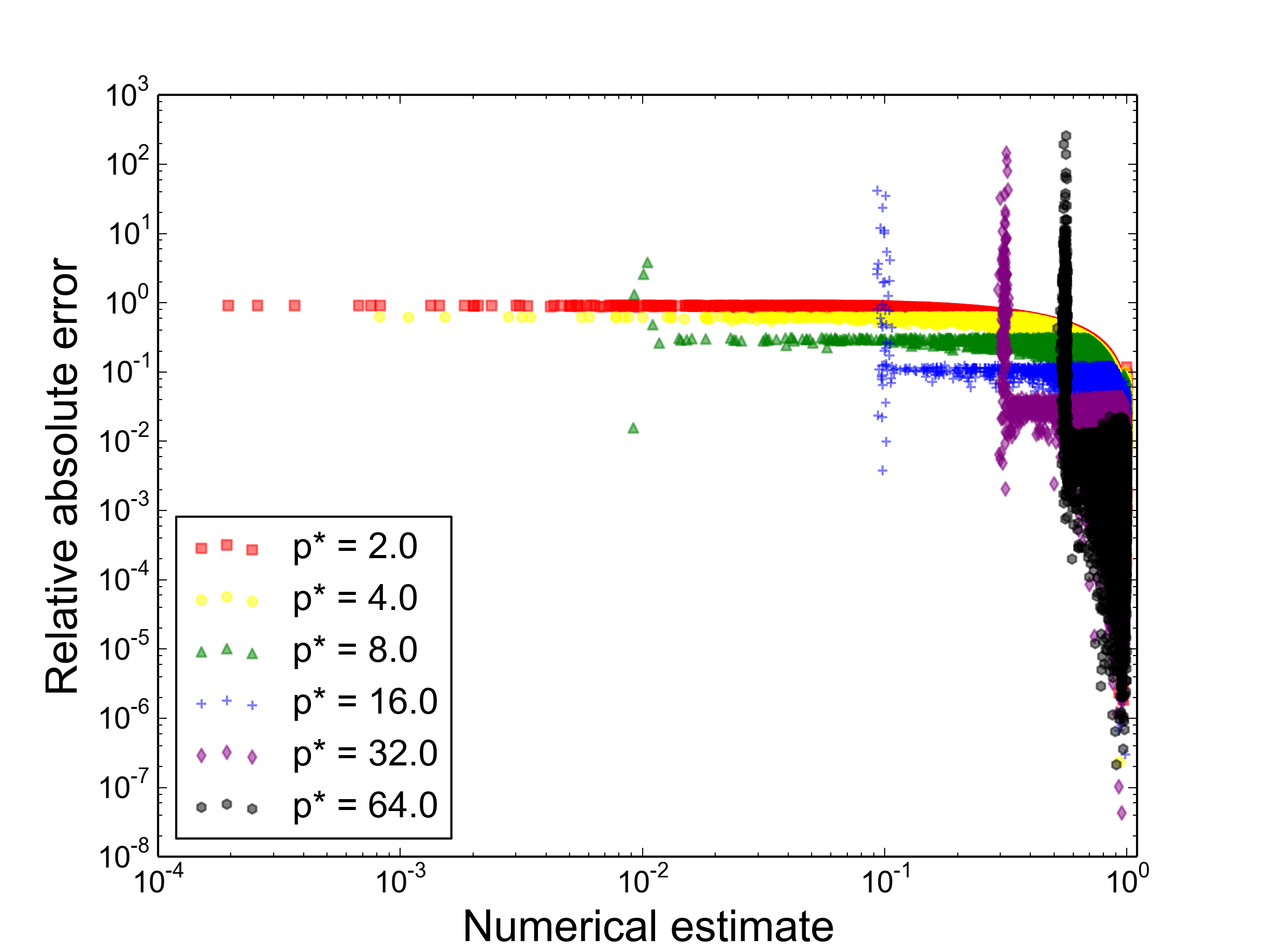}
  \end{tabular}

  \caption{{\bf The influence of the parameter $p^*$ on
      max-convolution accuracy.} For each $k \in \{256,512,1024\}$, 64
    replicate max-convolutions are performed to compare the relative
    error of the fast numerical method compared to the naive
    method. This is performed for different values of $p^*$; lower
    values of $p^*$ perform well when the exact result is close to
    zero, and higher values of $p^*$ perform better otherwise, and
    this relationship seems invariant of $k$ when the data are scaled
    in the manner presented (they are scaled so that the largest
    element has value 1).
  \label{figure:max-convolution-accuracy}}
\end{figure*}

A cursory empirical test of the numerical stability as a result of
$p^*$ and the vector length $k$ was performed. In {\bf
  figure~\ref{figure:max-convolution-accuracy}}, the numerical
stability was demonstrated on 64 random pairs of vectors for each
length $k \in 128, 256, 512, 1024$. For all $p^* \in \{ 2, 4, 8, 16,
32, 64\}$, and the relative absolute error of each element in the
result of the max-convolution is computed $|\frac{numerical[m] -
  exact[m]}{exact[m]}|$, where $numerical[m]$ and $exact[m]$ refer to
the value at index $m$ of the numerical and naive results
respectively. ({\bf Figure~\ref{figure:max-convolution-accuracy}})
demonstrates the relationship between $p^*$, $k$, the relative
absolute error and the magnitude of the exact result.

Qualitatively, optimizing the numerical performance involves
satisfying competing ideals: when $p^*$ is large, the problem solved
converges to the max-convolution inference, but underflow becomes
significant. When $p^*$ it's too small, non-maximal terms contribute
to the result (more similar to a standard convolution).

Although more sophisticated numerical analysis would almost certainly
yield larger improvements to the method, a simple improvement is
exploited: Generally the relative error is fairly low, but only
becomes high in these experiments when the exact value at that index
is close to zero (this is intuitive from the formula for relative
absolute error). Since underflow is the only numerical consideration
(because the values are normalized to the maximum), then it follows
that a result that is not close to zero at some index is convergent
when $p^*$ is substantially large (if it suffered from too much
underflow, then it would approach zero quickly). Therefore, when using
a high value of $p^*$, indices where the numerical solution is close
to zero indicate that the potential for numerical error, and suggest
that a smaller value of $p^*$ could be used for those indices. This
yields a further improvement where a more accurate result can be
constructed from two calls of {\bf
  algorithm~\ref{algorithm:max-convolution-revised}}; this improved
method is shown in {\bf
  algorithm~\ref{algorithm:max-convolution-piecewise}}, and runs in
roughly twice as many steps as {\bf
  algorithm~\ref{algorithm:max-convolution-revised}} (still $\in
O(k\log(k))$). Note that this piecewise method could be trivially
extended to use more than two values of $p^*$, increasing accuracy at
the expense of additional runtime (although, assuming the results
using the different values of $p^*$ are computed in decreasing order,
then the routine could potentially terminate once the result has been
estimated at all indices with adequate numeric stability).

\begin{algorithm}
  \caption{ {\bf Numerical max-convolution (piecewise normalized
      version)}, a numerical method to estimate the max-convolution of
    two PMFs or nonnegative vectors (further revised to strategically
    choose $p^*$). The parameters are two PMFs $\pmf_L$ and $\pmf_R$
    both with with $k$ outcomes (and both by definition
    nonnegative). The return value is a numerical estimate of the
    max-convolution $\pmf_L *_{\max} \pmf_R$. This algorithm calls the
    revised method $maxConvolutionRevised$ from {\bf
      algorithm~\ref{algorithm:max-convolution-revised}}. Note that
    the values of $p_{lower}^*$, $p_{higher}^*$, and $\tau$ specified
    here may be chosen to optimize performance on a specific
    application. }

  \label{algorithm:max-convolution-piecewise}
  \begin{small}
    \begin{algorithmic}[1]
      \Procedure{maxConvolutionPiecewise}{$\pmf_L$, $\pmf_R$}
      
      \State $p_{lower}^* \gets 4.0$
      \State $p_{higher}^* \gets 32.0$
      \State $\tau \gets 0.6$ \Comment{This is a two-way piecewise implementation; higher orders will be slower and more accurate}

      \State $f_{stable}[m] \gets maxConvolutionRevised(\pmf_L, \pmf_R, p_{lower}^*)$ \Comment{Numerically stable upper bound}
      \State $f_{aggressive}[m] \gets maxConvolutionRevised(\pmf_L, \pmf_R, p_{higher}^*)$ \Comment{More aggressive calculation}

      \For{$m$ to $k$, $m += 1$}
      \If{$f_{higher}[m] \geq \tau$}
      \State $f[m] \gets f_{aggressive}[m]$ \Comment{When high-$p^*$ result is large enough, underflow isn't an issue}
      \Else
      \State $f[m] \gets f_{stable}[m]$ \Comment{Where the high-$p^*$ result is small, use $f_{stable}$ to avoid underflow}
      \EndIf
      \EndFor

      \State \Return $f$
      \EndProcedure
    \end{algorithmic}
  \end{small}
\end{algorithm}

\subsection{Using fast numerical convolution to solve a probabilistic subset sum problem}

Lastly, a three-way piecewise implementation of max-convolution
similar to the one shown in {\bf
  algorithm~\ref{algorithm:max-convolution-piecewise}} but using $p^*
\in \{4,32,64\}$ (the Python code of this three-way piecewise method
is given in the accompanying Python demonstration code) is used to
solve a simulated probabilistic generalization of the the subset sum
problem. In this problem, $n=32$ people go shopping and each person $j
\ in \{1 , 2 , \ldots n\}$ buys exactly one of two items (the price of
the item they do purchase is $\mu_j^{(true)}$ and the price of the
item they do not purchase is $\mu_j^{(false)}$), where the costs of
both items for each person $j$ are unknown to us. Then, given fuzzy
knowledge about the costs of these items with all prices discretized
into $k=256$ bins (\emph{i.e.}, $\forall j,~ \pmf_{X_j}[\ell]$ where
$\ell \in \{0, 1, \ldots k-1\}$) and given fuzzy knowledge about the
total amount spent ($\pmf_M$, where $M = \sum_j X_j$), we try to infer
the amount spent by each person (\emph{i.e.}, for each person $j$,
estimating $\mu_j^{(true)}$).

Data are generated as follows: At each variable $X_j, j \in \{1, 2,
\ldots n\}$, two means are randomly sampled $\mu_j^{(true)},
\mu_j^{(false)} \sim uniform(0, k-1)$, and a discretized Gaussian PMF
is centered about each mean (the standard deviations of the Gaussians
are each sampled $\sigma_j^{(true)}, \sigma_j^{(false)} \sim
uniform(0, \frac{k}{10}) $). A vector proportional to the PMF
$\pmf_{X_j}[\ell]$ is computed using the sum of these Gaussians with a
vector of uniform noise $\alpha^{(in)}[\ell] \sim uniform(0,
0.0001)$. The likelihood distribution on the sum $M = \sum_j X_j$ is
generated by adding a Gaussian with mean $\sum_j \mu_j^{(true)}$ and
variance $0.005 \times (n k - (n-1))$ (\emph{i.e.}, $0.005\times$ the
possible number of outcomes for $M$), plus point-wise samples of
uniform noise $\alpha^{(out})[\ell] \sim uniform(0, 0.0001)$.

On each problem instance, likelihoods for all inputs $(\pmf_{Y_1},
\pmf_{Y_2}, \ldots \pmf_{Y_n})$ are computed twice, once using naive
max-convolution and once using the numerical method. Note that even
though these values of $n$ and $k$ do not appear particularly large,
the full max-convolution tree that they produce will will compute
several max-convolutions on the order of $O(k)$ and a few on the order
of $O(n\times k)$, which in this case can be $32\times 256 =
8192$. For this reason, computing the likelihood curve with the naive
result requires 159 seconds, while the fast numerical approach takes
0.935 seconds to compute a highly similar result.

A single likelihood distribution $\pmf_{Y_j}$ for one particular $j
\in \{1, 2, \ldots n\}$, which was computed using the naive method,
the fast numerical method is plotted in {\bf
  figure~\ref{figure:max-convolution-tree-accuracy}}. This figure also
demonstrates the utility of {\tt max-product} inference by also
showing the result from {\tt sum-product} inference, which is much
less informative (and does not have a mode close to the correct
answer.

\section{Discussion}

\begin{figure}
  \centering 
  \includegraphics[width=4in]{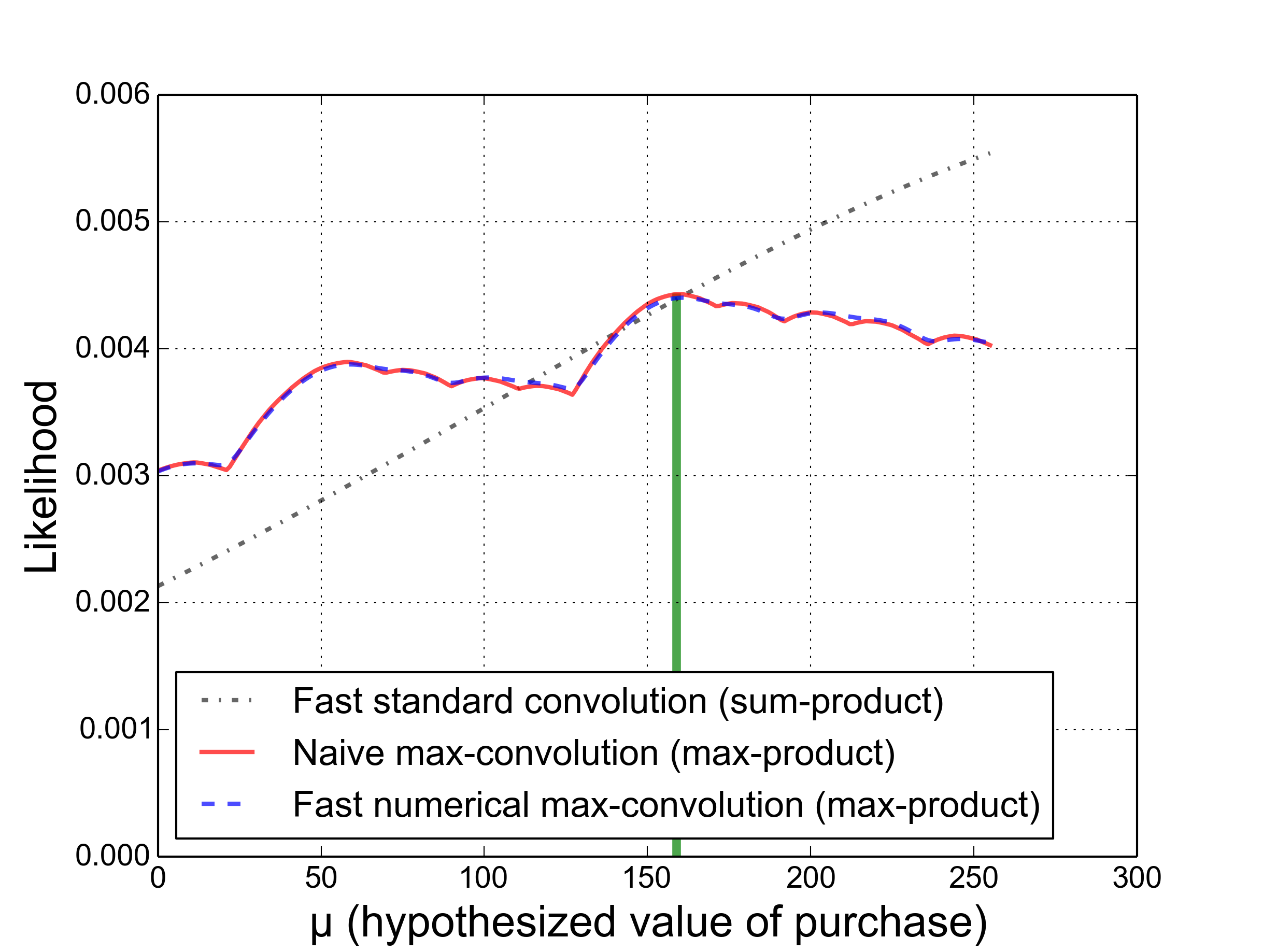}

  \caption{{\bf Using fast numerical max-convolution for {\tt
        max-product} inference.} A single likelihood distribution for
    one particular $j \in \{1, 2, \ldots n\}$ from the probabilistic
    generalization of the subset-sum problem is shown. This
    distribution is computed via a probabilistic max-convolution tree
    (the problem is solved twice, once with naive max-convolution and
    once via fast numerical max-convolution). The true mode value
    $\mu_j^{(true)}$ for that $j$ is indicated by the bar beneath the
    largest mode. To compare the results of different types of
    inference, the {\tt sum-product} result from the convolution tree
    (using the standard convolution operator) is also plotted; note
    that {\tt sum-product} inference produces a less discriminative
    likelihood curve because many possible joint events have diffused
    into it.
   \label{figure:max-convolution-tree-accuracy}}
\end{figure}

Although its ethos may ultimately limit the utility of this method to
numerical settings where small errors are tolerable (as opposed to the
to more general theoretical papers previously
mentioned\cite{bussieck:fast, bremner:necklaces}), numerical method
proposed here gives a simple and very fast estimate of the
max-convolution result, which could allow use of numerical
max-convolution (or {\tt max-product} inference on the sums and
differences between discrete distributions) in settings where it is
currently far too computationally expensive. The largest caveat to the
method is the inaccuracy that can occur due to numerical bottlenecks
(\emph{e.g.}, underflow); however, for many problems (\emph{e.g.},
practical applications of the {\tt max-product} inference problem in
{\bf figure~\ref{figure:max-convolution-tree-accuracy}}), the
numerical method is sufficient to perform high-quality inference, but
in a dramatically faster time.

Furthermore, the connection between the max-convolution problem and the
all-pairs shortest path problem from graph
theory\cite{bremner:necklaces} means that the fast numerical method
may be used to compute fast numerical approximations to that important
computer science problem. Such fast numerical estimates could
complement theoretical solutions to that problem.

A more in-depth theoretical analysis of the algorithm's error would
likely yield multiple opportunities to modify the algorithm in order
to decrease error. For example, one possible improvement could be
performed by using log-transformed real values: in log-transformed
space, raising to the power $p^*$ would be equivalent to scaling by
$p^*$, and would not produce significant underflow. Furthermore, the
operations required by FFT convolution could be performed in log-space
by translating the ring ($+$,$\times$) on real values to its
equivalent ring $(\log_+, \log_\times) = (\log_+, +)$ on
log-transformed values, where the operation $\log_+(x,y)$ is performed
by dividing out the greater of the two arguments $x$ (w.l.o.g.) and
then computing $\log(1+\frac{y}{x}) = \log(1+z)$ via Taylor series;
performing the Cooley-Tukey FFT on log-transformed values (or possibly
using a different FFT algorithm that is particularly well suited for
log-transformed values) could represent one route for improving the
accuracy.

In a similar vein, more sophisticated techniques for locally choosing
among a small number of values for $p^*$ (compared to a simple
piecewise function on two possible values of $p^*$). For instance,
under roughly uniform distributions of values in the two vectors being
max-convolved, the values closest to zero (and thus having higher
chance of having high relative absolute error) will occur more
commonly at the first and last indices of the numerical estimate
(because those indices take the maximum over smaller collections of
elements, and are thus more likely to be smaller values). Similar
attention could be put toward scaling the vectors prior to taking
elements to the $p^*$ (compared to the current procedure of dividing
by the maximum element value) may minimize the underflow on a large
number of points with large values. A most exciting possibility would
be that having an accurate estimate of the max-convolution result
somehow could be used to compute a more accurate result (\emph{e.g.},
using approximate results from different $p^*$ and exploiting the
property $\|\cdot\|_p \geq \|\cdot\|_{p+\delta}$ when $p \geq 1$ and
$\delta>0$). Such directions of future research could possibly solving
the max-convolution iteratively over a bounded or constant number of
subproblems where each subproblem requires $O(k\log(k))$, by first
computing initial estimates of the max-convolution result with the
numerical method presented here, and then using those initial
estimates to parameterize a subsequent call to the numerical method in
a manner reminiscent of the QR algorithm for
eigendecomposition\cite{francis1961qr, francis1962qr}. From {\bf
  algorithm~\ref{algorithm:max-convolution-piecewise}}, it seems
highly likely that there will be more ways by which an initial result
can be used to obtain a higher-accuracy result.

Furthermore, even pursuing methods for obtaining very high accuracy
with large values of $p^*$ may not always be of substantial
interest. Indeed, even if no improvement to accuracy is ever
presented, the design and parameterization of machine learning methods
(including graphical models) has traditionally been empirically
driven, and the use of exact {\tt max-product} inference is hardly
sacrosanct in every application (as opposed to inference somewhere
between {\tt sum-product} and {\tt max-product}). From this
perspective, rather than choosing $p^*$ as a static constant value
$p^*=1$ (\emph{i.e.}, performing {\tt sum-product} inference) or $p^*
\rightarrow \infty$ (\emph{i.e.}, performing {\tt max-product}
inference), $p^*$ could be viewed as a hyperparameter that is used to
position inference on a continuum somewhere between all joint events
contributing equally to the end result ({\tt sum-product}) and only
the best joint event contributing ({\tt max-product}), and
intermediate values of $p^*$ would establish a preference for the top
few joint events. In this sense, the value chosen for $p^*$ could be
driven by the data, and the problem of $p$-norm convolution (where a
finite $p$ is desired, rather than max-convolution where $p
\rightarrow \infty$) can already be solved with very high accuracy by
the proposed numerical method for any moderate choice of $p^*$.

\section{Availability}
\noindent A simple illustration of the fast numeric max-convolution
method in Python (using the {\tt numpy} package), including a
three-way piecewise implementation, is available at
\url{https://bitbucket.org/orserang/fast-numerical-max-convolution}.

\section*{Acknowledgements}
Thanks to Mattias Franberg and Ryan Emerson for the helpful comments.

\end{document}